\documentstyle[amscd,amssymb,verbatim,epsf,mathrsfs,graphicx,hyperref,pst-node,tikz-cd,fancyvrb,lipsum]{amsart}

\begin{document}
\makeatletter
\newcommand{\verbatimfont}[1]{\renewcommand{\verbatim@font}{\ttfamily#1}}
\makeatother
\newcommand{\R}{\mathbb{R}}
\newcommand{\C}{\mathbb{C}}
\newcommand{\om}{\omega}
\newcommand{\I}{{\mbox{I}}}
\newcommand{\II}{{\mbox{II}}}
\newcommand{\III}{{\mbox{III}}}
\newcommand{\Hess}{\mbox{\rm{Hess}}}
\newcommand{\Sing}{\mbox{\rm{Sing}}}
\newcommand{\Ind}{\mbox{\rm{Ind}}}
\theoremstyle{plain}
\newtheorem{Thm}{Theorem}
\newtheorem{Cor}{Corollary}
\newtheorem{Ex}{Example}
\newtheorem{Con}{Conjecture}
\newtheorem{Main}{Main Theorem}
\newtheorem{Lem}{Lemma}
\newtheorem{Prop}{Proposition}

\theoremstyle{definition}
\newtheorem{Def}{Definition}
\newtheorem{Note}{Note}
\newtheorem{Question}{Question}

\newtheorem{remark}{Remark}
\newtheorem{notation}{Notation}
\renewcommand{\thenotation}{}

\renewcommand{\rm}{\normalshape}%

\title{On Isolated Umbilic Points}
\author{Brendan Guilfoyle }
\address{School of STEM \\
          Munster Technological University \\
          Clash \\
          Tralee  \\
          Co. Kerry \\
          Ireland.}
\email{brendan.guilfoyle@@ittralee.ie}

\begin{abstract}
Counter-examples to the famous conjecture of Carath\'eodory, as well as the bound on umbilic index proposed by Hamburger, are constructed with respect to Riemannian metrics that are arbitrarily close to the flat metric on Euclidean 3-space. 

In particular, Riemannian metrics with a smooth strictly convex 2-sphere containing a single umbilic point are constructed explicitly, in contradiction with any direct extension of Carath\'eodory's conjecture to non-Euclidean metrics. Additionally, a non-Euclidean metric and embedded surfaces containing isolated umbilic points of any index are presented, violating Hamburger's umbilic index bound.

In both cases, it is shown that the metric can be made arbitrarily close to the flat metric. 
\end{abstract}

\maketitle


\vspace{0.1in}
\noindent {\em Dedicated to Karen Uhlenbeck on the occasion of her 75th birthday, in appreciation of her generous support during the years of cigar smuggling, bicycle riots and Zapatismo.}\footnote{A short video explaining the results of this paper can be found at the following link: \url{https://youtu.be/Wjja4PcMtxc}}
\vspace{0.1in}
\section{Isolated Umbilic Points}

Let $(M^3,g)$ be a smooth Riemannian 3-manifold and $S$ a smoothly embedded surface in $M$. The second fundamental form at a point, being a symmetric matrix, has two real eigenvalues, the principal curvatures. The {\it umbilic points}, where the principal curvatures are equal, are generically isolated, and away from them, the associated eigen-directions form a pair of orthogonal line fields on $S$.      

The eigen-directions determine the (unoriented) principal foliations on $S$, which have singularities at the umbilic points. For an isolated umbilic point $p\in S$ there is a well-defined half-integer index ${\mbox{Ind }}(p)\in{\mathbb Z}/2$ determined by the winding number of the principal foliation around the umbilic point.

In classical surface theory, where $(M^3,g)=({\mathbb R}^3,g_{0})$, $g_0$ being the flat Euclidean metric, surfaces which are umbilic at every point can be proven to be pieces of either spheres or planes \cite{Hopf_book}. At the other extreme, isolated umbilic points have many mysterious properties. 

Isolated umbilic points are the quarks of classical surface theory - they are contained on a string (a surface in ${\mathbb R}^3$), come in different flavours (the half-integer index) and are never seen alone in ${\mathbb R}^3$, at least on strictly convex surfaces (the Carath\'eodory Conjecture).

In addition, the distinction between smoothness and real analyticity of a surface can be intimately related to its isolated umbilic points, as has been observed from both a dynamic \cite{gak_hopf} and a stationary \cite{HandW1} \cite{Hopf} perspective. 

This paper demonstrates the failure of both the classical Carath\'eodory Conjecture and Hamburger's index bound when the metric is not flat, even for metrics arbitrarily close to Euclidean.

Recall that the infamous Carath\'eodory Conjecture asserts that a smooth strictly convex topological 2-sphere in Euclidean 3-space must have at least {\em two} umbilic points. Of course, if all of the umbilic points are isolated, the sum of the indices is the Euler characteristic of the 2-sphere $\chi(S)=2$, so there must be at least one umbilic point. If a 2-sphere had only one umbilic, its index would be 2 and so the Conjecture rules out spheres in which the principal foliation is that of, say, the standard dipole \cite{Laz}.

Hamburger proposed establishing the Carath\'eodory Conjecture by proving the stronger local bound ${\mbox{Ind }}(p)\leq1$ on the index of any isolated umbilic point \cite{Ham0}. In the case of real analytic surfaces, Hamburger established this bound \cite{Ham1} \cite{Ham2}. Subsequent published work sought unsuccessfully to simplify this proof \cite{klotz} \cite{titus}- for a recent contribution see \cite{ivanov}. Later, this proposed index bound was attributed to L\"{o}wner \cite{lowner}, but we will refer to it here as Hamburger's local index bound. 

Much subsequent research on the Conjecture has focused on this local index bound \cite{ando2} \cite{gutsanbring1} \cite{gutsanbring2} \cite{gutmercsanbring} \cite{gutasot} \cite{smythx1} \cite{smythandx2} \cite{sotmello} \cite{xavier}. Counter-examples which are not $C^2$-smooth have been constructed \cite{AndoFaU}, and related results discussed in  \cite{ghomihoward}.

One of the difficulties is that it has been unclear where exactly the problem lies: geometry, analysis or topology. The proof of the Conjecture in \cite{gak0} provides the answer: geometry. The proof hinges on the size of the Euclidean group and its action on the space of oriented geodesics. This is also behind the global to local argument that proves the index bound ${\mbox{Ind}}(p)<2$ \cite{gak00}, although this inequality has a topological nature.

In the above proof, the size of the isometry group implies that a certain related elliptic boundary value problem is Fredholm regular when there is only one umbilic point on a closed convex 2-sphere. This fact has been generalized to other similar elliptic boundary value problems in plane bundles, as long as there is a transitive isometric action \cite{gak000}. 

The final ingredient to the proofs of both global and local versions given in \cite{gak0} and \cite{gak00} is mean curvature flow with boundary in co-dimension two \cite{gak0000}. This is used to construct solutions of the elliptic boundary value problem and prove the existence of co-kernel, contradicting Fredholm regularity and establishing the Conjecture.

It is conceivable that the existence of holomorphic discs could be established by some PDE method other than mean curvature flow. In any event, the existence of co-kernel is sufficient to preclude a convex 2-sphere with a single umbilic point. It is also sufficient to prove a related conjecture of Toponogov on the umbilic points of complete planes \cite{gak_top} \cite{top}.

In the present paper we explore further this contention that it is geometry that lies at the heart of the Conjecture. From this view-point, there is no reason to believe that the Carath\'eodory Conjecture holds for a generic Riemannian 3-manifold, even ones close to the Euclidean metric, as they do not have sufficient isometries to fix a point. This would explain the subtlety and depth of the problem.

We demonstrate that this is indeed the case, by constructing counter-examples which show that neither the Carath\'eodory Conjecture nor the local index bound of Hamburger need hold for non-flat metrics, even ones that are arbitrarily close to the Euclidean metric $g_0$.

In detail, it is proven that isolated umbilic points of any index can be realized on surfaces in near-Euclidean metrics:
\vspace{0.1in}
\begin{Thm}\label{t:1}
For all $\epsilon>0$ and $k\in{\mathbb Z}/2$, there exists a smooth Riemannian metric $g$ on ${\mathbb R}^3$ and a smooth embedded surface $S\subset{\mathbb R}^3$ such that
\begin{itemize}
    \item[(A1)] $S$ has an isolated umbilic point of index $k$,
    \item[(B1)] $\|g-g_{0}\|^2\leq\epsilon$,
\end{itemize}
where $\|.\|$ is the $L_2$ norm on ${\mathbb R}^3$ with respect to the flat metric $g_0$.
\end{Thm}
\vspace{0.1in}

Turning to the global Carath\'eodory Conjecture it is proven that

\vspace{0.1in}
\begin{Thm}\label{t:2}
For all $\epsilon>0$, there exists a smooth Riemannian metric $g$ on ${\mathbb R}^3$ and a smooth strictly convex 2-sphere $S\subset{\mathbb R}^3$ such that
\begin{itemize}
    \item[(A2)] $S$ has a single umbilic point,
    \item[(B2)] $\|g-g_{0}\|^2\leq\epsilon$.
\end{itemize}
\end{Thm}
\vspace{0.1in}

The metrics $g$ are obtained by perturbing the Euclidean metric as follows. Fix a totally umbilic surface in flat ${\mathbb R}^3$ - a plane for Theorem \ref{t:1} and a round sphere for Theorem \ref{t:2}. The Euclidean metric is deformed, leaving the surface fixed, in such a manner that the induced metric on $S$ is preserved.

This means that the curvature introduced by the ambient metric, through the Gauss equation, deforms only the second fundamental form of $S$. The perturbations are sufficient to control the principal foliation of $S$ and can be scaled arbitrarily $C_0$ close to the flat metric. Bump functions with support in a collared neighbourhood of the surface can then be used to bring the metric $L_2$ close to $g_0$.

This paper is organized as follows. In the next section we give the geometric background for the computations that follow. Theorem \ref{t:1} is proven in Section \ref{s:3}, while Section \ref{s:4} contains the proof of Theorem \ref{t:2}.

\section{The Second Fundamental Form}\label{s:2}

Suppose that $S$ is a smooth oriented surface in a Riemannian 3-manifold $({\mathbb M}^3,g)$.  Choose an orthonormal frame $\{e_0,e_1,e_2\}$ along $S$ adapted to the surface in that $e_0$ is a unit normal, while $e_1$ and $e_2$ are unit tangent vectors to $S$. Clearly this is only well-defined up to a rotation about the normal, which, if we introduce complex null frames
\[
e_+={\textstyle{\frac{1}{\sqrt{2}}}}\left(e_1-ie_2\right) \qquad\qquad e_-={\textstyle{\frac{1}{\sqrt{2}}}}\left(e_1+ie_2\right),
\]
is given by $(e_+,e_-)\rightarrow (e^{i\theta}e_+,e^{-i\theta}e_-)$. 

At a point $p\in S$ the second fundamental form $\Pi :T_pS\times T_pS\rightarrow {\mathbb R}$ is defined
\[
\Pi(X,Y)=-g(\nabla_Xe_0,Y)=-g(\nabla_Ye_0,X),
\]
where $\nabla$ is the Levi-Civita connection of $g$.

The second fundamental form can be projected onto the basis $(e_+,e_-)$ which span $T_pS$, and its three real components become encoded in the quantities
\begin{equation}\label{e:spinco_def}
\sigma=\Pi(e_+,e_+) \qquad\qquad \rho=\Pi(e_+,e_-) ,
\end{equation}
where $\sigma$ is complex-valued and $\rho$ is real-valued.

The geometric significance of these for surface theory is the following:
\[
|\sigma|={\textstyle{\frac{1}{2}}}|\kappa_1-\kappa_2| \qquad\qquad \rho=-{\textstyle{\frac{1}{2}}}(\kappa_1+\kappa_2),
\]
where $\kappa_1,\kappa_2$ are the principal curvatures of $S$. The complex quantity $\sigma$ is called the {\it shear} or {\it astigmatism} and it vanishes at umbilic points. The real quantity $\rho$ is, up to a constant multiple, the mean curvature of $S$. 

The argument of $\sigma$ determines the principal directions, as measured against the complex null frame. In particular, if $S$ is given in local coordinates $(z,\bar{z})$ about the origin, and $\sigma=Hz^n\bar{z}^m$ for $n,m\in{\mathbb N}$ and $H(z,\bar{z})$ a real non-zero function, then the origin would be an isolated umbilic point of index $(m-n)/2$. 

\vspace{0.1in}
\section{Hamburger's local index bound}\label{s:3}

In this section we prove Theorem \ref{t:1} by showing that every index can be realized by umbilic points on surfaces in certain Riemannian metrics close to the Euclidean metric.

For all $\epsilon>0$ and each $k\in{\mathbb Z}/2$, a smooth Riemannian metric $g$ on ${\mathbb R}^3$ and smoothly embedded surfaces $S\subset{\mathbb R}^3$ will be constructed such that
\begin{itemize}
    \item[(A1)] $S$ has an isolated umbilic point of index $k$,
    \item[(B1)] $\|g-g_{0}\|^2\leq\epsilon$.
\end{itemize}

Endow ${\mathbb R}^3$ with standard coordinates $(x^1,x^2,x^3)$ and let $z=x^1+ix^2$. The Euclidean metric $g_0$ in coordinates $(z,\bar{z},t)$ is
\[
ds^2=dzd\bar{z}+dt^2.
\]
Consider the following metric $g$ defined 
\[
ds^2=dzd\bar{z}+dt^2+[\bar{\beta}(z,\bar{z})dz+\beta(z,\bar{z}) d\bar{z}]dt,
\]
where $\beta:{\mathbb R}^2\rightarrow{\mathbb C}$ is a smooth function. This has the following properties:
\vspace{0.1in}
\begin{Prop}\label{p:l2fl}
The metric $g$ is Riemannian iff $|\beta|^2<1$, while $\beta=0$ is the flat Euclidean metric. Moreover, the $L_2$-distance of $g$ to $g_0$ is
\[
\|g-g_{0}\|^2=2\|\beta\|^2=2\iiint_{{\mathbb R}^3}|\beta|^2d^3V.
\]
\end{Prop}
\begin{pf}
A straightforward computation shows that
\[
{\mbox{det}}g={\textstyle{\frac{1}{4}}}(1-\beta\bar{\beta}),
\]
which implies the first statement. 

The $L_2$ norm $\|.\|$ is defined with respect to the flat metric $g_0$, which means that 
\[
\|g-g_{0}\|^2=2\iiint_{{\mathbb R}^3}|g-g_{0}|^2d^3V,
\]
where, in coordinates with summation convention,
\[
|g-g_{0}|^2=(g_{ij}-g_{0ij})(g_{kl}-g_{0 kl})g_0^{ik}g_0^{jl}.
\]
The result then follows from verifying that
\begin{equation}\label{e:c0_fl}
|g-g_{0}|^2=(g_{ij}-g_{0 ij})(g_{kl}-g_{0 kl})g_0^{ik}g_0^{jl}=2|\beta|^2.
\end{equation}
\end{pf}
\vspace{0.1in}

Now consider the plane $S$ in ${\mathbb R}^3$ given by $t=0$. Clearly the metric $g$ restricted to this surface is flat since the induced metric is simply $dzd\bar{z}$. 

A natural null basis for the tangent space of $S$ is
\begin{equation}\label{e:e+fl}
e_+=\sqrt{2}\frac{\partial}{\partial z} \qquad\qquad e_-=\sqrt{2}\frac{\partial}{\partial \bar{z}},
\end{equation}
while the unit normal is easily computed to be
\begin{equation}\label{e:e0fl}
e_0=\frac{1}{\sqrt{1-\beta\bar{\beta}}}\left(\frac{\partial}{\partial t}-\beta\frac{\partial}{\partial z}-\bar{\beta}\frac{\partial}{\partial \bar{z}}\right).
\end{equation}
Computing the second fundamental form components on the null frame
\vspace{0.1in}
\begin{Prop}
For the surface $t=0$ we have
\begin{equation}\label{e:spinco_fl}
\sigma=\frac{2}{\sqrt{1-\beta\bar{\beta}}}\frac{\partial\bar{\beta}}{\partial z}
\qquad\qquad
\rho=\frac{1}{\sqrt{1-\beta\bar{\beta}}}\left(\frac{\partial\beta}{\partial z}+\frac{\partial \bar{\beta}}{\partial {\bar{z}}}\right).
\end{equation}
\end{Prop}
\begin{pf}
By definition 
\[
\sigma=\Pi(e_+,e_+)=-g(\nabla_{e_+}e_0,e_+)=-e_+^ie_+^l(\partial_ie_0^j+\Gamma_{ik}^je_0^k)g_{jl},
\]
\[
\rho=\Pi(e_+,e_-)=-g(\nabla_{e_+}e_0,e_-)=-e_+^ie_-^l(\partial_ie_0^j+\Gamma_{ik}^je_0^k)g_{jl}.
\]
Computing the Christoffel symbols $\Gamma{ij}^k$ of $g$ and using equations (\ref{e:e+fl})
 and (\ref{e:e0fl}) this reduces to the stated quantities.
 
 \end{pf}
\vspace{0.1in}
We are now in a position to construct metrics and surfaces with isolated umbilic points of any index. First fix $m,n\in{\mathbb N}$ and choose  $0<r_0<r_1<1$. For $\lambda>0$ define the function $\beta:{\mathbb R}^2\rightarrow{\mathbb C}$ by
\[
\beta(z,\bar{z})=\left\{ \begin{matrix}
\lambda z^n\bar{z}^m&\qquad\qquad {\mbox{for }}0\leq|z|\leq r_0\\
\lambda \Phi z^n\bar{z}^m&\qquad\qquad {\mbox{for }}r_0\leq|z|\leq r_1\\
0&\qquad\qquad {\mbox{for }}r_1\leq|z|
\end{matrix}\right.,
\]
where $\Phi$ is a smooth bump function with $0\leq\Phi\leq1$, $\Phi(|z|=r_0)=1$ and $\Phi(|z|=r_1)=0$.

This metric is Riemannian so long as $\lambda<1$ by Proposition \ref{p:l2fl} and the fact that
\begin{equation}\label{e:beta_fl}
|\beta|^2\leq \lambda^2.
\end{equation}

Now for the plane $t=0$ the second fundamental form component $\sigma$ in the neighbourhood of the origin is computed by substituting $\beta=\lambda z^n\bar{z}^m$ in the first of equations (\ref{e:spinco_fl}) yielding
\[
\sigma=\frac{2m\lambda z^{m-1}\bar{z}^n}{\sqrt{1-\lambda^2 |z|^{2m+2n}}}. 
\]
Thus the origin $z=0$ is an isolated umbilic point of index $k=(n-m+1)/2$. This establishes condition (A1) for suitable choice of  $n$ and $m$.

The construction so far has yielded a metric that can be made $C_0$ close to the flat metric at each point by choosing $\lambda$ small enough, since by equations (\ref{e:c0_fl}) and (\ref{e:beta_fl}) 
\[
|g-g_{0}|^2=2|\beta|^2\leq2\lambda^2.
\]
To extend this to an $L_2$ estimate, introduce a further bump function $\Psi:{\mathbb R}\rightarrow{\mathbb R}$ in the metric so that
\[
ds^2=dzd\bar{z}+dt^2+\Psi(t)[\bar{\beta}(z,\bar{z})dz+\beta(z,\bar{z}) d\bar{z}]dt,
\]
where $\beta:{\mathbb R}^2\rightarrow{\mathbb C}$ is the smooth function above and
\[
\Psi(t)=\left\{ \begin{matrix}
1&\qquad\qquad {\mbox{for }}|t|\leq \epsilon/4\\
f(t)&\qquad\qquad {\mbox{for }}\epsilon/4\leq|t|\leq \epsilon/2\\
0&\qquad\qquad {\mbox{for }} \epsilon/2\leq|t|
\end{matrix}\right.,
\]
$f(t)$ being a smooth bump function with $0\leq f(t)\leq1$. The $L_2$ difference between the flat metric and this bumped metric is
\[
\|g-g_{0}\|^2=2\iiint_{{\mathbb R}^3}f(t)^2|\beta|^2d^3V\leq2\int_{-\epsilon/2}^{\epsilon/2}dt\iint_{|z|^2\leq1}\lambda^2dzd\bar{z}
\leq2\lambda^2\epsilon\pi.
\]

Thus choosing $\lambda<1/\sqrt{2\pi}$ ensures condition (B1) holds.

\vspace{0.1in}
\section{Convex spheres with a single umbilic point}\label{s:4}

In this section we prove Theorem \ref{t:2} by showing that there exist near-Euclidean Riemannian metrics containing smooth strictly convex topological 2-spheres with a single umbilic point. 

That is, for all $\epsilon>0$, there exists a smooth Riemannian metric $g$ on ${\mathbb R}^3$ and a smooth strictly convex 2-sphere $S\subset{\mathbb R}^3$ such that
\begin{itemize}
    \item[(A2)] $S$ has a single umbilic point,
    \item[(B2)] $\|g-g_{0}\|^2\leq\epsilon$.
\end{itemize}
thus showing that Caratheodory's Conjecture need not hold for non-flat metrics.

Consider ${\mathbb R}^3$ with standard spherical polar coordinates $(R,\theta,\phi)$ about the north pole. Introduce the complex coordinate $\xi=\tan(\theta/2)e^{i\phi}$ on the sphere.

The Euclidean metric $g_0$ in coordinates $(R,\xi,\bar{\xi})$ is
\[
ds^2=dR^2+\frac{2R^2}{(1+\xi\bar{\xi})^2}d\xi d\bar{\xi}.
\]
Consider the following metric $g$
\[
ds^2=dR^2+\frac{2R^2}{(1+\xi\bar{\xi})^2}d\xi d\bar{\xi}+\frac{R}{1+\xi\bar{\xi}}[\bar{\beta}(\xi,\bar{\xi})d\xi+\beta(\xi,\bar{\xi}) d\bar{\xi}]dR,
\]
where $\beta:{S}^2\rightarrow{\mathbb C}$. 
To extend this over the whole 2-sphere, change coordinates to the antipodal coordinates on the 2-sphere by
\[
\xi'=-\frac{1}{\bar{\xi}},
\]
and follow through the metric coordinate change. One finds that the metric retains the same form with
\[
\beta'(\xi',\bar{\xi}')=\frac{{\xi}'}{(1+\xi'\bar{\xi}')\bar{\xi}'}\bar{\beta}\left(-\frac{1}{\bar{\xi}'},-\frac{1}{{\xi}'}\right).
\]
This metric has the following properties:
\vspace{0.1in}
\begin{Prop}\label{p:l2sph}
The metric $g$ is Riemannian iff $|\beta|^2<1$, while $\beta=0$ is the flat Euclidean metric. 
The $L_2$-distance of $g$ to $g_0$ is
\[
\|g-g_{0}\|^2=2\|\beta\|^2=2\iiint_{{\mathbb R}^3}|\beta|^2d^3V.
\]
\end{Prop}
\begin{pf}
A computation shows that
\[
{\mbox{det }}g=\frac{4R^4(1-\beta\bar{\beta})}{(1+\xi\bar{\xi})^4},
\]
which implies the first statement. 

The proof of the second statement  is similar to the proof of the analogous statement in Proposition \ref{p:l2fl}.
\end{pf}
\vspace{0.1in}

Now consider the sphere $S$ in ${\mathbb R}^3$ given by $R=R_0>0$. Clearly the metric $g$ restricted to this surface is round. Thus the natural orthonormal basis for the tangent space is
\[
e_+=\frac{1+\xi\bar{\xi}}{\sqrt{2}R_0}\frac{\partial}{\partial \xi} \qquad\qquad e_-=\frac{1+\xi\bar{\xi}}{\sqrt{2}R_0}\frac{\partial}{\partial \bar{\xi}},
\]
while the unit normal is computed to be
\[
e_0=\frac{1}{\sqrt{1-\beta\bar{\beta}}}\left[\frac{\partial}{\partial R}-\frac{1+\xi\bar{\xi}}{2R_0}\left(\beta\frac{\partial}{\partial \xi}+\bar{\beta}\frac{\partial}{\partial \bar{\xi}}\right)\right].
\]
Computing the second fundamental form components on the null frame
\vspace{0.1in}
\begin{Prop}
For the surface $R=R_0$ the second fundamental form components on the null frame is
\begin{equation}\label{e:spinco_sp}
\sigma=\frac{1}{R_0\sqrt{1-\beta\bar{\beta}}}\frac{\partial}{\partial \xi} \left[(1+\xi\bar{\xi})\bar{\beta}\right],
\end{equation}
\begin{equation}
\rho=-\frac{1}{R_0\sqrt{1-\beta\bar{\beta}}}\left[4-(1+\xi\bar{\xi})^2\left[\frac{\partial}{\partial\xi}\left(\frac{\beta}{1+\xi\bar{\xi}}\right)+\frac{\partial}{\partial{\bar{\xi}}}\left(\frac{\bar{\beta}}{1+\xi\bar{\xi}}\right)\right]\right].
\end{equation}
\end{Prop}
\vspace{0.1in}
We are now in a position to construct the metric and 2-sphere with a single umbilic point. For $\lambda>0$ define the function $\beta:{S}^2\rightarrow{\mathbb C}$ by
\[
\beta(\xi,\bar{\xi})=\frac{\lambda\bar{\xi}}{(1+\xi\bar{\xi})^2}.
\]
It is easy to prove that
\[
|\beta|^2\leq \frac{3^3}{4^4}\lambda^2,
\]
so by Proposition \ref{p:l2sph} the metric is Riemannian so long as
\[
\lambda^2\leq \frac{4^4}{3^3}.
\]

For the 2-sphere $R=R_0$, the second fundamental form components are
\[
\sigma= \frac{\lambda}{R_0\sqrt{(1+\xi\bar{\xi})^4-\lambda^2\xi\bar{\xi}}}
\qquad\qquad
\rho= \frac{2(1+\xi\bar{\xi})^2+3\lambda\xi\bar{\xi}\cos2\theta}{2R_0\sqrt{(1+\xi\bar{\xi})^4-\lambda^2\xi\bar{\xi}}},
\]
where $\xi=|\xi|e^{i\theta}$.

There are no umbilic points for finite values of $\xi$. Thus the only umbilic point is at the south pole $|\xi|\rightarrow\infty$. To understand the point at the south pole, change coordinates to the antipodal coordinates on the sphere as above and compute
\[
\sigma'= \frac{\lambda\bar{\xi}^{'4}}{\sqrt{(1+\xi'\bar{\xi}')^4-\lambda^2\xi'\bar{\xi}'}},
\]
which gives a smooth isolated umbilic at $\xi'=0$ of index $2$. 

To see that $S$ is strictly convex, consider the determinant of the second fundamental form
\[
\kappa=\rho^2-\sigma\bar{\sigma}=\frac{[2(1+\xi\bar{\xi})^2+\lambda(1+3\xi\bar{\xi}\cos2\theta)][2(1+\xi\bar{\xi})^2-\lambda(1-3\xi\bar{\xi}\cos2\theta)]}{R_0[(1+\xi\bar{\xi})^4-\lambda^2\xi\bar{\xi}],
},
\]
and note that it is positive for $\lambda<1$. This establishes condition (A2).

The construction so far has yielded a metric that can be made $C_0$ close to the flat metric at each point by choosing $\lambda$ small enough, since
\[
|g-g_{0}|^2=2|\beta|^2\leq2\frac{3^3}{4^4}\lambda^2.
\]
To extend this to an $L_2$ estimate, introduce a further bump function $\Psi:{\mathbb R}\rightarrow{\mathbb R}$ in the metric so that
\[
ds^2=dR^2+\frac{2R^2}{(1+\xi\bar{\xi})^2}d\xi d\bar{\xi}+\frac{R\Psi(R)}{1+\xi\bar{\xi}}(\bar{\beta}(\xi,\bar{\xi})d\xi+\beta(\xi,\bar{\xi}) d\bar{\xi})dR,
\]
where $\Psi$ is similar to the function earlier, but now radially and centred at $R=R_0$:
\[
\Psi(R)=\left\{ \begin{matrix}
1&\qquad\qquad {\mbox{for }}|R-R_0|\leq \epsilon/4\\
f(R)&\qquad\qquad {\mbox{for }}\epsilon/4\leq|R-R_0|\leq \epsilon/2\\
0&\qquad\qquad {\mbox{for }} \epsilon/2\leq|R-R_0|
\end{matrix}\right.,
\]
\vspace{0.1in}
where again $f(R)$ is a smooth bump function with $0\leq f(R)\leq1$. 

The $L_2$ difference between the flat metric and this bumped metric is
\begin{align}
\|g-g_{0}\|^2&=2\iiint_{{\mathbb R}^3}\Psi(R)^2|\beta(\xi,\bar{\xi})|^2d^3V\nonumber \\
&\leq2\frac{3^3}{4^4}\lambda^2\int_{R_0-\epsilon/2}^{R_0+\epsilon/2}R^2dR\iint_{S^2}\frac{2d\xi d\bar{\xi}}{(1+\xi\bar{\xi})^2}\nonumber \\
&\leq2\frac{3^3\epsilon^2\pi}{4^4}\lambda^2(12R_0^2+\epsilon^2).\nonumber
\end{align}
Thus choosing 
\[
\lambda<\frac{\sqrt{2}^7}{\sqrt{3^3\pi(12R_0^2+\epsilon^2)}},
\]
ensures condition (B2) holds.

\end{document}